\documentclass{amsart}
\usepackage{amssymb, amscd}
\numberwithin{equation}{section}

\def\CC{{\mathbb C}}

\def\PP{{\mathbb P}}
\def\QQ{{\mathbb Q}}

\def\ZZ{{\mathbb Z}}

\def\Hcal{{\mathcal H}}

\def\Mcal{{\mathcal M}}

\def\Ocal{{\mathcal O}}

\def\Xcal{{\mathcal X}}

\newcommand\proofsquare{\nobreak\hfill \hbox{%
\vrule height 5pt
\kern-.4pt
 \vbox{%
\hrule width 5pt depth0pt height.4pt
 \kern4.6pt \hrule  }
\kern-3.75pt
\vrule height 5pt}\kern1pt
\par}

\newtheorem{theorem}{Theorem}[section]

\newtheorem{proposition}[theorem]{Proposition}
\newtheorem{corollary}[theorem]{Corollary}

\newtheorem{definition-lemma}[theorem]{Definition-Lemma}

\theoremstyle{definition}

\theoremstyle{remark}

\begin{document}

\title[On Tautological Classes on Moduli of K3]{
Note on Tautological Classes \\
of Moduli of K3 Surfaces
\footnote{\tt k3classes.tex\today}}
\author{Gerard van der Geer}
\address{Korteweg-de Vries Instituut, Universiteit van
Amsterdam, Plantage Muidergracht 24, 1018 TV Amsterdam, The Netherlands}
\email{geer@science.uva.nl}
\author{Toshiyuki  Katsura}
\address{Graduate School of Mathematical Sciences\\
 University of Tokyo \\ Komaba, Meguro-ku\\
 Tokyo\\
153-8914 Japan}
\email{tkatsura@ms.u-tokyo.ac.jp}


\begin{abstract}
In this note we prove some cycle class relations on moduli
of K3 surfaces.
\end{abstract}

\maketitle

\begin{section}{Introduction}
\label{sec: intro}
This note deals with a few properties of tautological classes
on moduli spaces of K3 surfaces. Let ${\Mcal}_{2d}$ denote a
moduli stack of K3 surfaces over an algebraically closed
field with a polarization of degree $2d$ prime to the characteristic
of the field.
The Chern classes of the relative cotangent bundle 
$\Omega^1_{{\Xcal}/{\Mcal}}$
of the universal K3 surface ${\Xcal}_{2d}={\Xcal}$ 
define classes $t_1$ and $t_2$ in the Chow groups
${CH}^i_{\QQ}({\Xcal}_{2d})$ of the universal K3 surface over ${\Mcal}_{2d}$.
The class $t_1$ is the pull back from ${\Mcal}_{2d}$ of the first Chern
class $v=c_1(V)$ of the Hodge line bundle $V=
\pi_*(\Omega_{{\Xcal}/{\Mcal}}^2)$. We use Grothendieck-Riemann-Roch
to determine the push forwards of the powers of $t_2$. These are 
powers of $v$. We then prove that $v^{18}=0$ in the Chow group with rational
coefficients of ${\Mcal}_{2d}$ 
We show that this implies that a complete subvariety of ${\Mcal}_{2d}$
has dimension at most $17$ and that this bound is sharp.
These results are in line with those for moduli of abelian varieties.
There the top Chern class $\lambda_g$ of the Hodge bundle vanishes in
the Chow group with rational coefficients. The idea is that
if the boundary of the Baily-Borel compactification has co-dimension $r$
then some tautological class of co-dimension $r$ vanishes. Our result
means that $v^{18}$ is a torsion class. It would be very interesting 
to determine the order of this class as well as explicit representations
of this class as a cycle on the boundary, cf., \cite{E-vdG1,E-vdG2}.

\begin{section}{The Moduli Space ${\Mcal}_{2d}$}
Let $k$ be an algebraically closed field. We consider the moduli
space ${\Mcal}_{2d}$ of polarized K3 surfaces over $k$ with a primitive
polarization of degree $2d$.
This is a 19-dimensional algebraic space. Over the complex numbers
we can describe it as an orbifold quotient 
$\Gamma_{2d}\backslash \Omega_{2d}$,
where $\Omega_{2d}$ is a bounded symmetric domain and 
$\Gamma_{2d}$ an arithmetic
subgroup of $SO(3,19)$ obtained 
as follows. Consider the lattice 
$U^3\oplus E_8^2$, where $U$ denotes a hyperbolic plane
and $E_8$ the usual rank $8$ lattice.
Let $h$ be an element of this lattice with $\langle h, h \rangle = 2d$.
Then $L_{2d}=h^{\bot} \cong U^2 \bot E_8^2 \bot \ZZ u$ with
$\langle u , u \rangle = -2d$ is a lattice of signature $(2,19)$ and we put
$$\Omega_{2d}=\{ [\omega]  \in \PP(L_{2d} \otimes \CC) :
\langle \omega , \omega \rangle =0, \, \langle \omega , \bar{\omega} \rangle >0 \}.
$$
The group $\Gamma_{2d}$ is the automorphism group of $L_{2d}$.
It acts on $\Omega_{2d}$ and the quotient (an orbifold) 
is the analytic space of  ${\Mcal}_{2d}$.
It is well-known by Baily-Borel
that the sections of a sufficiently high power of 
$V$ give an embedding of $\Gamma_{2d} \backslash \Omega_{2d}$ 
as a quasi-projective variety.

\end{section}

\end{section}
\bigskip
\begin{section}{GRR applied to the sheaf $\Omega_{\Xcal /M}^i$}
In order to determine the push forward $\pi_*(t_2^a)$ under
$\pi: {\Xcal}_{2d} \to {\Mcal}:={\Mcal}_{2d}$
we apply Grothendieck-Riemann-Roch to the structure sheaf 
of the universal (polarized) K3-surface $\pi:{\Xcal}\to {\Mcal}$. 
We work in the Chow ring with rational coefficients. 
We have
$$
{\rm ch}(\pi_!{\Ocal}_{\Xcal})= \pi_*({\rm ch}({\Ocal}_{\Xcal}) {\rm
Td}^{\vee}(\Omega_{\Xcal / \Mcal }^1))
= \pi_*({\rm Td}^{\vee}(\Omega^1_{{\Xcal}/ {\Mcal}})).
$$
As to the left-hand-side we have
$\pi_!{\Ocal}_{{\Xcal}}= 1 + V^{\vee}$,
where $V=R^0\pi_*\Omega^2_{\Xcal / \Mcal}$ is the line bundle 
with fibre $H^0(X, \Omega_X^2)$ over $[X]$. We write
$v$ for the first Chern class of this  bundle on ${\Mcal}$. So the
left hand side is
$ 1 +e^{-v}$.
For the right hand side, remark that $\Omega^1_{{\Xcal}/ {\Mcal}}$
has as determinant a line bundle that is trivial on each K3-surface
that is a fibre of $\pi$. Therefore, this line bundle is a pull back
from ${\Mcal}$ and we can identify it with $\pi^*(V)$. If
we denote by $t_i=c_i(\Omega_{\Xcal /\Mcal}^1)$
the Chern classes of $\Omega_{\Xcal /\Mcal}^1$, the
right hand side has the following form
$$
\pi_*(1-t_1/2 +  (t_1^2+t_2)/12 - t_1t_2/24 + ...).
$$
Comparing the terms of degree $0$ gives
$1+1= 24/12$ since $c_1^2(X)=0$ and $c_2(X)=24$ for a K3 surface.
For the terms of degree $1$ we find:
$-v = \pi_*(-t_1t_2)/24 =  -v\cdot (\pi_*(t_2)/24)$
and this is in agreement. Degree $2$ terms yield
$\pi_*(t_2^2)= 88\, t_1^2$.
This checks with the next term:
$$
v^3/6= {1 \over 1440}\pi_*(3t_2^2t_1-t_2t_1^3).
$$
Continuing this way we can determine $\pi_*(t_2^j)$ for all $j\geq 1$.
More precisely, put $B(x)=x/(1-e^{-x})$ and write $\gamma_1, \gamma_2$
for the Chern roots of $\Omega_{\Xcal /\Mcal}^1$. Then
$$
{\rm Td}^{\vee}(\Omega_{\Xcal /\Mcal}^1)=B(\gamma_1)B(\gamma_2)=
 \sum_{n, j \colon 0\leq 2j \leq n} c(n,j)(\gamma_1+\gamma_2)^{n-2j}
(\gamma_1 \gamma_2)^j
$$
with $t_1=\gamma_1+\gamma_2$ and $t_2=\gamma_1 \gamma_2$
and the Riemann-Roch identity says that if $\pi_*(t_2^{n+1})=a_nv^{2n}$
the $a_n$ satisfy the relation:
$$
\sum_{j\geq 0} a_j \, c_{n,j} =
\left\{
\begin{array}{l}
1/(n-2)! \quad \mbox{for}~ n\equiv 0 (\bmod \, 2), n>2,\\
2\quad \mbox{for}~ n=2.\\
\end{array} \right.
$$
Denoting $\pi_*(t_2^{n+1})=a_n v^{2n}$
we find the following values for $a_n$.
$$
\begin{matrix}
n & \pi_*(t_2^{n+1})/v^{2n} \cr
\hline
0 & 24  \cr
1 & 88  \cr
2 & 184 \cr
3 & 352 \cr
4 & 736 \cr
5 & 1295488/ 691  \cr
6 & 4292224 / 691   \cr
7 & 68418650624/ 2499347 \cr
8 & 17412311922527744/ 109638854849\cr
9 & 22654813560476770158592/ 19144150084038739\cr
\end{matrix}
$$
\smallskip
\noindent
\begin{proposition} Write $\pi_*(t_2^{n+1})=a_n v^{2n}$ for $n\geq 0$.
The generating function
$$
A(t)= \sum_{n=0}^{\infty}
 a_n\, 	t^n = 24 + 88\, t+184\,  t^2+ \ldots
$$
is characterised uniquely by the property that for every $n>0$ the
coefficient of $t^{2n-1}$ in
$$
{2-t\over 1-t}A\left( {-t^2\over 1-t} \right)
$$ 
is equal to $4n/B_{2n}$.  Here $B_{m}$ denotes the $m$-th Bernouilli number.
\end{proposition}
Although the numbers $a_n$ are defined for all $n \geq 0$ they have
a geometric interpretation for $n\leq 9$ only, apparently.

\bigskip

We now apply Grothendieck-Riemann-Roch to the sheaf 
$\Omega^1_{\Xcal/\Mcal}$, or
equivalently, to its dual $\Theta_{\Xcal / \Mcal}$. It says
$$
{\rm ch}(\pi_{!}\Theta_{\Xcal /\Mcal})= 
\pi_*({\rm ch}(\Theta_{\Xcal /\Mcal}) {\rm Td}^{\vee}(\Omega^1_{\Xcal /\Mcal})).
$$
Note that in the left hand side $\pi_!\Theta_{\Xcal / \Mcal}=
R^1\pi_*\Theta_{\Xcal /\Mcal}$
since a K3 surface has no non-zero vector fields, \cite{R-S}.
Since the push forward of powers of $t_2$ are powers of $v$ and $t_1=\pi^*(v)$
we see that ${\rm ch}(\pi_{!}\Theta_{\Xcal /\Mcal})$ is a polynomial in $v$. 
This can be determined by looking at cohomology once 
we show that the tautological ring of ${\Mcal}_{2d}$ is $\QQ[v]/v^{18}$.

Note that the fibre of $R^1\pi_* \Theta_{\Xcal / \Mcal}$ is $H^1(X,\Theta_X)$,
the space of infinitesimal deformations of $X$. The tangent space to $M$ at
$[X]$ can be identified with the orthogonal complement of $h$, the hyperplane
class in $H^1(X,\Omega_X^1)=H^1(X,\Theta_X)$. 
On the other hand we know that
by Hodge theory the following description for this tangent space. Let
$$
0 \subset F^2 \subset F^1 =(F^2)^{\bot } \subset H^2_{dR}
$$
be the Hodge filtration on $H^2_{dR}(X)$ and let $h$ be the hyperplane class
which gives a section of $H^2_{\rm dR} \otimes O_{\Mcal}$. Then the tangent space to $\Mcal$
can be identified with ${\rm Hom}(F^2,(F^1\cap h^{\bot})/F^2)$. Using
the cup product we can identify $(F^2)^{\vee}$ with $H_{\rm dR}^2/F^1$,
i.e., in the Grothendieck group we have $[H^2_{\rm dR}]=V^{\vee}+ V^{\bot}$,
where we identify $F^2$ with $V$. We now restrict to the
orthogonal subbundle $h^{\bot}$ of the hyperplane class $h$ whose
class in the Grothendieck group is $[H^2_{\rm dR}]-1$.
Therefore we find $[\Theta_{\Mcal}]=[(H^2_{\rm dR}-1-V-V^{-1})\otimes V^{-1}]$.

\begin{proposition} In the Grothendieck group of ${\Mcal}$ we have the
relation $[\Theta_{\Mcal}]=[H^2_{\rm dR}-1] \otimes V^{-1}-1-V^{-2}$.
\end{proposition}

In view of the Gauss-Manin connection on $H_{\rm dR}^2$ we see
that the Chern classes vanish in cohomology
and that the total Chern class
of the bundle $F^1$ on $\Mcal$ is $1/(1-v)$ and
${\rm ch}(\Theta_{\Mcal})=-1+21e^{-v}-e^{-2v}$ 
and in particular we find $c_1(\Theta_{\Xcal/ \Mcal})=-19v$. We already saw that
the total Chern class of $R^1\pi_*(\Omega^1_{\Xcal / \Mcal})$ is
$1/(1-v^2)$ and so its first Chern class vanishes. This checks with
global duality $(R^1\pi_*\Theta_{\Xcal / \Mcal})^{\vee}
\cong (R^1\pi_*\Omega_{\Xcal/ \Mcal}^1)\otimes V$.
\end{section}

\begin{section}{Vanishing of tautological classes in characteristic zero}
Let $II_{3,19}$ be the unique even unimodular lattice of signature $(3,19)$
and let $S$ be some Lorentzian sublattice of $II_{3,19}$, say of signature
$(1,m)$. 
Recall that an $S$-K3 surface $X$ is a K3 surface with a fixed primitive 
embedding of $S$ into the Picard group such that the image of $S$ 
contains a semi-ample class, i.e., a class $D$ such that $D^2>0$ and
$D\cdot C\geq 0$ for all curves $C$ on the K3 surface $X$, cf., \cite{BKPS}.
The period space $Y$ of marked $S$-K3 surfaces is an orbifold which is a
quotient of a hermitean symmetric domain of dimension $19-m$
by an arithmetic subgroup of ${\rm Aut}(S^{\bot})$.

\begin{theorem} For $m\leq 16$ the cycle class $v^{18-m}$ vanishes in
the Chow group $CH^{18-m}_{\QQ}(Y)$ with rational coefficients.
\end{theorem}

\begin{proof}
By imposing a level structure we can replace our period space by a 
finite cover and assume that we are working with a fine moduli space.

The proof is by descending induction on $m$. 
For $m=16$ the period domain can be identified with the Siegel upper half space 
${\Hcal}_2$ and the orbifold $Y$ can be viewed as a moduli space of 
abelian surfaces. It thus carries a natural vector bundle, the Hodge bundle 
$\pi_*(\Omega^1_{{\Xcal}/Y})$ with Chern classes $\lambda_1$ and $\lambda_2$.
One shows that $\lambda_1=v$ by comparing the factors of automorphy
or by noticing that $H^0(X,\Omega_X^2)\cong
\wedge^2 H^0(\Omega_X^1)$ for an abelian surface. 
Furthermore, it is known that $\lambda_1^2$
vanishes by \cite{vdG}, Prop.\ 2.2. 
We conclude that $v^2$ vanishes. 

The induction step is now provided by Theorem 1.2 of \cite{BKPS}.
There exists a modular form $\Phi$ of weight $k\geq 12$ whose 
zero-divisor is of the form $\sum m_i W_i$ with $m_i \in \ZZ_{>0}$
with orbifolds $W_i$ that are images in $Y$ of
quotients $\Gamma_{L_t}\backslash \Omega_{L_t}$ under finite maps.
Here $\Omega_{L_t}$ is a hermitean symmetric domain of dimension one
less than the original domain 
$\Omega$ and the quotient parametrizes a family of $S'$-K3-surfaces 
with $S'\supset S$ of signature $(1,m+1)$. By induction we know
that on each of the orbifolds $\Gamma_{L_t}\backslash \Omega_{L_t}$ 
the class $v^{17-m}$ vanishes. 
The zero-divisor of $\Phi$ represents the class
$k \, v$. We thus find that a non-zero multiple of 
$v^{18-m}=v\cdot v^{17-m}$ vanishes. 
\end{proof}
In characteristic $0$ we can use the existence
of the Satake compactification whose boundary is $1$-dimensional 
to conclude
that intersecting twice with a sufficiently general
hyperplane yields a complete $17$-dimensional
subvariety of ${\Mcal}$. 
Since by Baily-Borel the class $v$ is ample this shows
that $v^{17}\neq 0$.
\begin{corollary}
The tautological ring of ${\Mcal}_{2d}$ is $\QQ[v]/(v^{18})$.
\end{corollary}
                                                                                
\begin{corollary}
The maximal dimension of a complete subvariety of ${\Mcal}_{2d}$ is $17$.
\end{corollary}

In positive characteristic the locus of K3-surfaces with
height $\geq 3$ defines a complete subvariety of
dimension $17$, cf.\ 
\cite{vdG-K}.

If ${\Mcal}^*$ is the Baily-Borel compactification of ${\Mcal}$
then the `boundary' is a $1$-dimensional cycle. In the Chow group
$CH^{18}_{\QQ}({\Mcal}^*)$ the class $v^{18}$ is represented by
a $1$-cycle with support on the boundary. 

\end{section}

 \end{document}